\documentclass{elsarticle}

\def\LaTeX{L\kern-.36em\raise.3ex\hbox{a}\kern-.15em
    T\kern-.1667em\lower.7ex\hbox{E}\kern-.125emX}

 \usepackage{latexsym,alltt,fleqn}
\usepackage{url}
\usepackage{textcomp}
 \usepackage{graphics}
 \usepackage{epsfig}
 \usepackage{times}
 \usepackage{amsmath}
  \usepackage{amssymb}
\usepackage[colorlinks]{hyperref}
 \usepackage[ps,curve,2cell]{xypic}
 \UseAllTwocells
 \usepackage[active]{srcltx}
 \usepackage{xy}
 \xyoption{v2}
 \usepackage{multirow}
\def\Hlang{\textit{H}}
\def\Hspec{\textit{Hspec}}
\def\Hver{\textit{Hver}}

\newcommand{\op}{\mathrm{op}}

\def\HI{\mathcal{HI}}

\def\CAT{\mathbf{C\!A\!T}}

\newcommand{\Nom}{\mathrm{Nom}}

\newcommand{\Mod}{\mathit{Mod}}
\newcommand{\Dom}{\mathit{Dom}}
\newcommand{\SEN}{\mathit{Sen}^{\HI}}
\newcommand{\Sen}{\mathit{Sen}}

\newcommand{\Sign}{\mathit{Sign}}

\def\p{\mathrm{th}}

\def\Set{\mathbf{Set}}

\def\FOL{\mathcal{F\!O\!L}}

\newcommand{\w}{w}

\def\implies{\mathbin{\Rightarrow}}

\newcommand{\D}{\mathcal{D}}
\newcommand{\I}{\mathcal{I}}
\newcommand{\co}{\,\colon\;}
\newcommand{\ra}{\rightarrow}

\newtheorem{example}{Example}[section]

\newcommand{\smallitems}{\setlength{\parsep}{-0.5ex plus 0.5ex}%
                        \setlength{\itemsep}{-0.5ex plus 0.5ex}%
                        \vspace{-1ex}
                        }
\DeclareMathAlphabet{\mathbb}{U}{msb}{m}{n}
\DeclareSymbolFont{ams}{U}{msa}{m}{n}
\DeclareSymbolFontAlphabet{\mathams}{ams}
\DeclareMathSymbol{\filter}{\mathams}{ams}{22}

\usepackage{xcolor}

\begin{document}

\title{Introducing H, 
  an institution-based formal specification and
  verification language} 

\author{R\u{a}zvan Diaconescu}
\ead{Razvan.Diaconescu@imar.ro}
\address{Simion Stoilow Institute of Mathematics of the Romanian
Academy, Bucharest, Romania}

\journal{} 

\begin{abstract}
This is a short survey on the development of the formal specification
and verification language \Hlang\ with emphasis on the scientific
part.
\Hlang\ is a modern highly expressive language solidly based upon
advanced mathematical theories such as the internalisation of
\emph{Kripke semantics} within \emph{institution theory}. 
\end{abstract}

\maketitle


\section{Introduction}

\Hlang\ is a language for formal specification and verification that has
emerged out of a theoretical effort spread over a period of more than
25 years.
It has been provisionally implemented as a prototype running system
during 2017 -- 2018 \cite{codescu:hybrid-2019} and following this
implementation a number of succesful case studies have been reported 
\cite{codescu:hybrid-2019,Tutu-Chirita-Lopes-Fiadeiro:LSBSS-2019}. 
\Hlang\ is designed as a two-component system: 
\begin{itemize}

\item \Hspec\ -- the specification language. 
  This is an institution-based language in the sense that is
  parameterised over a variety of base logic systems captured as
  institutions (in the sense of the \emph{institution theory} of
  Goguen and Burstall \cite{ins}). 
  The role of the base logics refer to the specification of the data
  part of the system. 
  For the dynamics part, \Hspec\ employs the essentials of modal
  logic.
  This design is based upon the understanding that the essential
  ingredients of modal logic (both at the syntactic and at the
  semantic level) are independent of the base level logic.  
  The fact that \Hspec\ is parameterised by base logics gives it
  unparalleled specification expressivity and power since the most
  appropriate logic for the data part may be chosen. 
  Moreover the list of base logics is open, any new logic may be
  added when convenient. 

\item \Hver\ -- a collection of verification tools and methods for
  specifications developed with \Hspec. 
  At this moment \Hver\ contains only one such method which is based
  on translations to first order logic. 
  However in the future we plan to extend \Hver\ with other methods
  and tools. 

\end{itemize}
A web page of \Hlang\ is currently maintained at
\begin{quotation}
\href{http://imar.ro/~diacon/forver/forver.html}
{\url{http://imar.ro/~diacon/forver/forver.html}}.  
\end{quotation} 
The \emph{\Hlang\ concept}, the core of \Hlang\ represented by its
scientific foundation, is based on developing an abstract Kripke semantics within
the institution theory of Goguen and Burstall \cite{ins}.  
This is the core of \Hlang, from which the vision of \Hlang\ had
emerged, its development being the result of a sustained mathematical
effort reported in a series of papers from which the most
representative are \cite{ks,QVHybrid,EncHybrid}.  
In this way \Hlang\ may be a typical example of a formal method that has
emerged out of initially purely theoretically motivated work. 

In this paper we survey the development of \Hlang, with emphasis on its
most important aspects. 
This goes as follows:
\begin{enumerate}

\item In the first part, which is also the main part of the survey, we
  present the mathematical foundations of \Hlang\ and its basic design. 
  As mentioned above this consists essentially of the internalisation
  of Kripke semantics in abstract institutions.
  This concept has been developed very gradually over many years and
  in our paper we will survey the most important moments of this
  development. 
  I believe that presenting the \Hlang\ concept in this way has several
  benefits for the reader, including the possibility to understand the
  flow of ideas behind the \Hlang\ concept.
  This may be quite difficult if we chose to present only the end
  result of this rather complex process. 

\item The next section is dedicated to the current implementation of
  \Hlang. This is only a temporary implementation that is based on the Hets
  system \cite{mossakowski-maeder-lutich2007}. In long term, our
  vision for \Hlang\ is to have an independent implementation. 

\item In the last section we discuss very briefly some case studies
  that have been formally specified and verified with \Hlang.

\end{enumerate}

The readership must be familiar with some very basic category theory
concepts, which are now quite commonly used in some areas of formal
methods. 
Some familiarity with institution theory and modal logic may be also
quite helpful.  

\section{The \Hlang\ concept}

The broader scientific context of the \Hlang\ concept is the theory of
institutions of Goguen and Burstall \cite{ins}.
The narrower context is the development of Kripke semantics within
abstract institutions \cite{ks,QVHybrid,EncHybrid}. 
In this section we first give a very brief presentation of institution
theory, and then survey the development process of the
institution-theoretic Kripke semantics. 

\subsection{Institutions}

The model theory oriented axiomatic approach by Goguen and Burstall
to the notion of a logical system \cite{ins} that is based on the
notion of \emph{institution} has started a line of important
developments of adequately abstract and general approaches to the
foundations of software specifications and formal system development
(see \cite{sannella-tarlecki-book}) as well as a modern version of
very abstract model theory (see \cite{iimt}).   
One of the main original motivations for introducing institution
theory was to respond to the explosion in the population of logics in
use in computing almost four decades ago, a situation that continues
today perhaps at an accelerated pace. 
These days the concept of institution lies at the foundations of several
formal modern specification languages and environments such as Maude
\cite{maude-book}, CASL \cite{casl2002} or CafeOBJ \cite{caferep},
Hets \cite{mossakowski-maeder-lutich2007} etc. 
In the area of formal specification and verification the contribution
of the institution-theoretic approach to modularity and heterogeneity
are priceless.  
Let us recall the notorious concept of institution: 

\noindent\fbox{%
    \parbox{\textwidth}{%
An  \emph{institution} $\I = 
\big(\Sign^{\I}, \Sen^{\I}, \Mod^{\I}, \models^{\I}\big)$ consists of 
\begin{itemize}\smallitems

\item a category $\Sign^{\I}$ whose objects are called
  \emph{signatures},

\item a sentence functor $\Sen^{\I} \co \Sign^{\I} \ra \Set$
  defining for each signature a set whose elements are called
  \emph{sentences} over that signature and defining for each signature
  morphism a \emph{sentence translation} function, 

\item a model functor $\Mod^{\I} \co (\Sign^{\I})^{\op} \ra \CAT$
  defining for each signature $\Sigma$ the category
  $\Mod^{\I}(\Sigma)$ of \emph{$\Sigma$-models} and $\Sigma$-model
  homomorphisms, and for each signature morphism $\varphi$ the
  \emph{reduct} functor $\Mod^{\I}(\varphi)$,  

\item for every signature $\Sigma$, a binary 
  \emph{$\Sigma$-satisfaction relation}
  $\models^{\I}_{\Sigma} \subseteq |\Mod^{\I} (\Sigma)|
  \times \Sen^{\I} (\Sigma)$, 

\end{itemize}
such that for each morphism 
$\varphi\co\Sigma \rightarrow \Sigma' \in \Sign^{\I}$, 
the \emph{Satisfaction Condition}
\begin{equation}
M'\models^{\I}_{\Sigma'} \Sen^{\I}(\varphi)(\rho) \text{ if and only if  }
\Mod^{\I}(\varphi)(M') \models^{\I}_\Sigma \rho
\end{equation}
holds for each $M'\in |\Mod^{\I} (\Sigma')|$ and $\rho \in \Sen^{\I} (\Sigma)$.
\[
\xymatrix{
    \Sigma \ar[d]_{\varphi} & \big|\Mod^{\I}(\Sigma)\big|
    \ar@{-}[r]^-{\models^{\I}_{\Sigma}} & 
    \Sen^{\I}(\Sigma) \ar[d]^{\Sen^{\I}(\varphi)} \\
    \Sigma' & \big| \Mod^{\I}(\Sigma')\big| \ar[u]^{\Mod^{\I}(\varphi)} 
    \ar@{-}[r]_-{\models^{\I}_{\Sigma'}} & \Sen^{\I}(\Sigma')
  }
\] 

    }%
}

\

The literature (e.g. \cite{iimt,sannella-tarlecki-book}) shows myriads
of logical systems from computing or from mathematical logic captured
as institutions. In fact, an informal thesis underlying institution
theory is that any `logic' may be captured by the above
definition. While this should be taken with a grain of salt, it
certainly applies to any logical system based on satisfaction between
sentences and models of any kind.
In concrete institutions, typically the signatures are structured
collections of symbols, the sentences are inductively defined from
atoms by using sentence building operators, the sentence translations
(along signature morphisms) rename symbols, the models interpret the
symbols of the signatures as sets and functions, the reducts
``forget'' interpretations of some symbols, and the satisfaction is
defined inductively on the structure of the sentences in Tarski's style
\cite{tarski44}. 

Here we refrain from presenting examples of logical systems captured
as institutions since the institution theory literature abounds of
such examples.
Instead let us just point out that the process of defining particular
logical systems as institutions is not necessarily a trivial one since
one may have to reconsider and give a serious fresh thought to
concepts such as signature morphisms, variables, quantifiers, etc.  
This rethinking of various concepts may have to do very much with the
intended applications, such as formal specification.
For example, from the specification perspective the concept of
signature morphism has to be much more general than what is usually
employed in conventional logic, in order for the mathematics to work
the variables require a kind of qualifications that are inspired from
the practice of specification languages, etc. 
Some of these issues have been discussed in extenso in
\cite{UniLogCompSci}.

\subsection{Kripke semantics in institutions}

The semantics for modal logics, known as \emph{Kripke semantics} was
introduced in \cite{kripke59}. 
The origin of the development of Kripke semantics in institutions --
often refereed to as `modalization of institutions' -- lies in some
research undertaken within the group of the late Professor Joseph Goguen
at Oxford in the early nineties regarding institutions for modal
logics. 
First there was the realisation of the fact that the model
amalgamation properties in modal logic institutions are a direct
consequence of the respective properties in the base logics, such as
propositional or first order logic.  
From there it followed the idea that each modal logic institution has
an underlying simpler base institution and that the Kripke models may
be defined uniformly on the basis of the models in a base institution.  
However it took over a decade to see the first paper on this published
\cite{ks}, mainly due to a rather complicated refereeing process. 

In \cite{ks} -- which may be considered the seminal paper for the
\Hlang\ concept -- we have introduced the first version of Kripke
semantics in abstract institutions first by considering a ``base''
institution $\I$ and then by building a ``modal'' institution $\HI$ on
top of $\I$.
This construction has several components:
\begin{enumerate}

\item An extension of the syntax of $\I$. 
While the signatures stay the same, new sentences are built from the
sentences of $\I$ by iteration of sentences building operators such as
the usual Boolean operators, quantifiers, and modalities. 

\item Kripke models built from the models of $\I$. 

\item The definition of a modal satisfaction relation between the
  Kripke models and the new sentences.

\end{enumerate}
Now let us review these three components of the construction of $\HI$
one by one.

\subsubsection{The syntax of $\HI$}

For any signature $\Sigma$, the set $\Sen^{\HI} (\Sigma)$ of the
$\Sigma$-sentences of the ``modal'' institution $\HI$ is the least set
closed under the following operations: 

\

\noindent\fbox{%
    \parbox{\textwidth}{%
\begin{itemize}

\item $\Sen^{\I} (\Sigma)\subseteq \SEN(\Sigma)$;

\item $\rho \star \rho' \in \SEN(\Sigma)$ for any
  $\rho, \rho'\in \SEN(\Sigma)$ and any
  $\star \in \{\vee,   \wedge,\implies\}$,  

\item  $\neg \rho \in \SEN(\Sigma)$, for any $\rho\in \SEN(\Sigma)$,

\item $[\lambda](\rho_1,\dots,\rho_n), \langle \lambda \rangle(\rho_1,\dots,\rho_n)
\in \SEN(\Sigma)$, \\
for any $\lambda\in \Lambda_{n+1}, \rho_i \in \SEN(\Sigma)$, 
$i\in \{1,\dots, n\}$;

\item $(\forall \chi)\rho, (\exists \chi)\rho\in \SEN(\Sigma)$, 
for any $\rho \in \SEN(\Sigma')$ and 
$\chi\co\Sigma \rightarrow \Sigma'\in \D$;


\end{itemize}
    }%
}

\

Some explanations are necessary:
\begin{itemize}

\item The first condition says that each sentence of the base
  institution becomes automatically a sentence of the
  ``modal''institution. 

\item The second and the third conditions does the Boolean connectors
  on the sentences of the ``modal'' institution. 
  Note that the sentences of the base institution may also involve
  Boolean connectors, in this case it is important to distinguish
  between the Boolean connectors at the base level and at the modal
  level since in general their effects may be differ. 

\item The fourth condition introduces modalities as sentences building
  operators. Here $\Lambda_n$ means the set of modalities of arity
  $n$, which may be thought just as relation symbols. 
  For now the modalities are \emph{not} considered part of the
  signatures, they are rather fixed. 

\item Quantifiers are considered in the institution-theoretic manner,
  via designated signature morphisms (see for example \cite{iimt} for
  details). Conventional concrete quantifiers would correspond to
  those signature morphisms that are in fact extensions of signatures
  with variables. 
  So, not each signature morphisms may be used in quantifiers, those
  that are designated for such use form a so-called
  \emph{quantification space} which is the $\mathcal{D}$ from the last
  condition above. 
  This concept represents an axiomatic approach to quantifiers that
  considers coherence properties with respect to translations along
  signature morphisms; it has been defined first in \cite{cond} and
  given this name in \cite{QVHybrid}. 

\item Like with the Boolean connectors we have to carefully
  distinguish between quantifiers  the level of $\HI$ and those that
  come with the sentences of the base institution as their effects may
  differ. 

\item The general institution-theoretic feature of the quantifiers,
  namely that they support higher-order quantification (up to what the
  concrete concept of signature supports) applies also here. 
  So, depending on how we chose $\mathcal{D}$ we may have first order,
  or second order, or even higher order quantifiers. 

\end{itemize}

\subsubsection{Kripke models}

The models of $\HI$ are Kripke models defined on the basis of the
models of the base institution $\I$:

\

\noindent\fbox{%
    \parbox{\textwidth}{%
Given a signature $\Sigma$, a \emph{Kripke model} $(M,W)$ consists of 
\begin{itemize}

\item a set $|W|$  -- called the set of the ``possible worlds''; 

\item for each $\lambda \in \Lambda_{n}$, a relation $W_\lambda
  \subseteq |W|^n$; and 

\item a mapping $M \co |W| \to |\Mod^{\I} (\Sigma)|$. 

\end{itemize}
    }%
}

\

\noindent
So, for each $w\in W$, $M_w$ is a model of the base institution $\I$. 
Moreover, $W = (|W|, (W_\lambda)_{\lambda\in \Lambda})$ is called the
\emph{Kripke frame} of $(M,W)$. 

However in order for the quantifications to work properly, usually the
models $M_w$ have to share something. 
For example in the concrete case of first order modal logic it is
quite common to require that the first order logic models $M_w$ that
are part of a Kripke model share their underlying sets and the
interpretations of the variables. 
At the level of the abstract institutions this condition has been
expressed in a general way in \cite{ks} as 

\

\noindent\fbox{%
    \parbox{\textwidth}{%
$\beta_\Sigma (M_w) = \beta_\Sigma (M_v)$ for all $w,v \in |W|$.
    }%
}

\

\noindent
where $\beta_\Sigma \co \Mod(\Sigma) \to \Dom(\Sigma)$ is a functor
satisfying some rather mild technical conditions (we omit them here). 

\subsubsection{The ``modal'' satisfaction}

The satisfaction relation that relate the syntax of the ``modal
institution'' $\HI$ to its semantics is defined by following the usual
institution theoretic definitions and in Tarski's style by recursion
on the structure of the sentences, the recursion base being the
satisfaction in the base institution $\I$.
For each Kripke model $(M,W)$ and each $w \in |W|$ we define a
``local'' satisfaction relation as follows:

\

\noindent\fbox{%
    \parbox{\textwidth}{%
\begin{itemize}
		
\item $(M,W)\models^w \rho$ iff $M_w\models^{\I} \rho$; 
when $\rho \in \Sen^{\I}(\Sigma)$,

\item $(M,W)\models^w \rho \wedge \rho' $ iff  $(M,W)\models^\w \rho$  and  
$(M,W)\models^\w \rho'$, and similarly for the other Boolean
connectors; 

\item $(M,W)\models^w [\lambda](\xi_1,\dots,\xi_n)$ iff for any 
$(\w,\w_1,\dots,\w_n) \in W_\lambda$ we have that
$(M,W)\models^{w_i}\rho_i$ for some $1\leq i\leq n$.

\item $(M,W)\models^w 
    \langle\lambda\rangle(\xi_1,\dots,\xi_n)$ iff there exists 
    $(\w,\w_1,\dots,\w_n)
    \in W_\lambda$ such that and $(M,W)\models^{w_i}\xi_i$ for any  
    $1\leq i\leq n$.

\item $(M,W)\models^w (\forall \chi) \rho$ iff 
$(M',W') \models^w \rho$ for any $(M',W')$ such that 
$\Mod(\chi)(M',W') = (M,W)$,   

\item $(M,W)\models^w (\exists \chi) \rho$ iff 
$(M',W') \models^w \rho$ for some 
$(M',W')$ such that $\Mod(\chi)(M',W') = (M,W)$, and 

\end{itemize}
    }%
}

\

\noindent
Under these definition in \cite{ks} it has been proved that $\HI$ is an
institution where the satisfaction $(M,W) \models \rho$ is defined on
the basis of the ``local'' satisfaction by $(M,W) \models^w \rho$ for
all $w\in |W|$.  
Moreover, the adequacy of this construction has been tested against
some deep model theoretic results including a very general
``modal'' ultraproducts theorem and its compactness consequences. 
Although in \cite{ks} we have not used multi-modalities (i.e. the
relations from $\Lambda$) but instead used the more familiar $\Box$
and $\Diamond$, this difference is insignificant, being just a matter
of form. 

Note that $\HI$ in fact represents a class of institutions rather than
a single institution because of the several parameters involved in its
construction.
Besides the base institution $\I$ of course, we also have the
modalities $\Lambda$, the quantification space $\D$ and the sharing
functor $\beta$. 
From this perspective a notation such as $\HI (\Lambda,\D,\beta)$
appears as more appropriate, however this is rather heavy so we
usually stick to the simpler version when the involved parameters are
clear.  

The usual modal logic institutions arise immediately as examples of
$\HI$. For instance modal propositional logic arises when considering
propositional logic as base institution (eventually stripped of the
Boolean connectors) and with $\D$ being trivial, while first order
modal logic arises when considering atomic first order logic as the
base institution (i.e. first order logic stripped off the Boolean
connectors and off the quantifiers) and $\D$ consisting of the
extensions of the signatures with appropriate variables.  
However the potential of the construction of $\HI$ goes much beyond
that of known examples of modal logics because it frees modal logic
from its conventional base. 
For example, at the base level it is possible to have partial functions
with various kinds of sharing (an interesting one from an
\cite{EncHybrid} would consider the sharing only of the definition
domains of the partial functions). 
A more intriguing example is given by the possibility to iterate
this construction for a number of times, obtaining hierarchical modal
logics. 

The construction of the ``modal institution'' of \cite{ks} is quite
emblematic for all other developments in the area and constitutes the
very basis for \Hspec\ as both the syntax and the semantics of \Hspec\
are based on this construction, but subject to some further additions
that will be presented below. 
 
\subsection{Adding nominals}

An important development in the area of institution-theoretic Kripke
semantics is the extension of the theory of \cite{ks} with the
ingredients of the so-called ``hybrid logic''.
Hybrid logics \cite{blackburn2000} are a brand of modal logics that 
provides appropriate syntax for the Kripke semantics in a 
simple and very natural way through the so-called \emph{nominals}. 
Historically, hybrid logic was introduced in \cite{prior} and further
developed in works such as 
\cite{Passy:1991:ECD:116976.116979,Areces01bringingthem,torben} etc.
The name ``hybrid logics'' was coined by Blackburn, but I consider
this an uninspired choice leading to confusions because of at least
two reasons.  
On the one hand this name does not suggest in any way the reality,
namely that ``hybrid logics'' is a sub-brand of the modal logics. 
In fact the difference between the two is rather minor because
technically it boils down only to a simple syntactic addition, whilst
they share the same semantics.   
On the other hand the term ``hybrid'' has a clear meaning in ordinary
language, which is difficult to relate to the corresponding brand of
logics. 
In spite of all these considerations, the terminology ``hybrid
logics'' is already established in the literature, and even the
name \Hlang\ owes to it.  

The presence of nominals brings in several advantages from the point
of view of formal specification and verification, such as the
possibility of explicit reference to specific states of the model and
a better more uniform proof theory. 
All these specification benefits have called for an extension of the
original theory of \cite{ks}, a first attempt in this direction being
\cite{HybridIns}.  
That had been technically a rather straightforward enterprise, which
is briefly presented below.   

\subsubsection{Upgrading the signatures} 

At the level of the signatures of $\HI$ we add the nominals, so after
this addition a signature consists of a pair $(\Nom, \Sigma)$ where
$\Nom$ is a set (of nominal symbols) and $\Sigma$ is a signature of
the base institution $\I$.  
This had been a good moment to include also the modalities $\Lambda$
in the signatures, a move that is specification oriented.  
When specifying dynamics of systems it is necessary to have user
defined modalities. 
Therefore a signature in $\HI$ is now a triple
$(\Nom,\Lambda,\Sigma)$. 

\subsubsection{Upgrading the sentences}

The collection of the sentence building operators gets expanded with: 

\

\noindent\fbox{%
    \parbox{\textwidth}{%
\begin{itemize}

\item $\Nom \subseteq \SEN(\Nom,\Lambda,\Sigma)$;

\item $@_i \rho\in \SEN(\Nom,\Lambda,\Sigma) $ for any $\rho \in
  \SEN(\Nom,\Lambda,\Sigma)$ and $i \in \Nom$; 

\end{itemize}
    }%
}

\

\noindent
Then there is the issue of upgrading the quantification by allowing
quantifications over the nominals.
For this we have to consider $\D$ as a quantification space for the
upgraded signatures, but one which does not have any effect on the
modalities. 
Thus the quantification building operators get upgraded to:

\

\noindent\fbox{%
    \parbox{\textwidth}{%
\begin{itemize}

\item $(\forall \chi)\rho, (\exists \chi)\rho\in \SEN(\Nom,\Lambda,\Sigma)$, 
for any $\rho \in \SEN(\Nom',\Lambda,\Sigma')$ and 
$\chi\co (\Nom,\Lambda,\Sigma) \rightarrow (\Nom',\Lambda,\Sigma') \in \D$;

\end{itemize}
    }%
}

\subsubsection{Upgrading the semantics}

The upgrade of the concept of Kripke models is very simple, just
interpret the new syntactic entities by extending $W$ with
interpretations for the nominals. 
So for each $i \in \Nom$ we have a designated element $W_i \in |W|$.   

\subsubsection{Upgrading the satisfaction relation}

This upgrade adds the semantics of the new building operators as
follows:

\

\noindent\fbox{%
    \parbox{\textwidth}{%
\begin{itemize}
		
\item $(M,W)\models^w i$ iff  $W_i=w$; when $i\in \Nom$,
		
\item $(M,W)\models^w @_j \rho$ iff $(M,W)\models^{W_j} \rho$. 

\end{itemize}
    }%
}

\

\noindent
We can see that the upgrade of the construction of $\HI$ from
\cite{ks} in the direction of nominals is technically very
straightforward. 
This is one of the reasons the paper \cite{HybridIns} may be
considered as only a minor contribution to the general development of
the \Hlang\ concept.
But there are other more serious reasons for this evaluation.
Due to being a conference paper -- and therefore being quite rushed
and suffering from severe space limitations -- the authors had to
scrap some crucial features of the original construction from
\cite{ks}, such as the sharing at the level of the Kripke models. 
One of the dramatic consequences of this simplification -- called the
``free hybridisation'' -- was that the quantification became
nonfunctional, thus reducing a lot the specification power of the
formalism.    
However these shortcomings have been corrected in the journal paper
\cite{QVHybrid}, which may be considered as the first paper
addressing the extension of \cite{ks} with nominals in a proper way. 

\subsection{More general constraints}

In \cite{QVHybrid} an ultimate very general axiomatic approach to the 
\emph{constraints} on Kripke models had been proposed. 
This approach captures a wide variety of constraints, such as various
sharing constraints or constraints on the shape of the Kripke frames
(such as reflexivity, transitivity, etc.).  
It is for instance more general and more accommodating than the
sharing constraints defined in \cite{ks}.  
Let us recall from \cite{QVHybrid}:

\

\noindent\fbox{%
    \parbox{\textwidth}{%
A \emph{constrained $\HI$-model functor} is a sub-functor 
$\Mod^C \subseteq \Mod^{\HI}$ such that it reflects weak amalgamation
for the designated pushout squares corresponding to the quantification
space $\D^{\I}$ at the level of the base institution $\I$ (that is
obtained by ``forgetting'' the nominals part from $\D^{\HI}$).
The models in $\Mod^C$ are called \emph{constrained $\HI$-models}.
    }%
}

\

\noindent 
We omit here detailed explanations concerning the technical condition
on weak amalgamation as the interested reader may consult
\cite{QVHybrid} or \cite{EncHybrid}. 
Informally, the meaning of the reflection condition of 
is that in the case of the designated pushout squares used in
quantifications the amalgamation of constrained models yields a 
constrained model. 
The role of this condition, which is rather mild in the applications,
is to ensure that the constrained models support smoothly the
quantifications. 
At the end we get a `sub-institution' of $\HI$ with constrained Kripke
models that is denoted $\HI^C$. 

\subsection{\Hspec\ and $\HI^C$}

The definition of \Hspec\ sticks closely to the construction of $\HI^C$,
being just a realisation of this construction as a specification
language. 
The following ideas underlie the definition of \Hspec:
\begin{itemize}

\item The syntax of \Hspec\ comes on two layers. 
The ``upper'' layer follows the definition of the signatures and of
the sentences of $\HI^C$, which become \Hspec\ declarations.
The ``lower'' layer follows the definition of the signatures and of
the sentences of the base institution $\I$, which is the most
important parameter the respective specification. 
In principle there is almost absolute freedom about the ``lower''
layer, in practice however we have to commit to something concrete,
usually to something that already exists in the realm of current
specification languages. 
For example CASL \cite{casl2002} may be used in many situations
because its underlying institution is a rather complex one which
includes many logical features, such as Boolean connectors, first
order quantifications, partial functions, etc. 

\item The semantics of a \Hspec\ specification is the class of the
  constrained $\HI^C$ models (Kripke models) that satisfy the axioms
  declared in the respective specification. 

\item In \Hspec, currently the constraints on the Kripke models are
  specified in two ways.  
  Either by ``rigidity'' declarations for the syntactic entities
  (sorts, operations, relations, etc.) that are meant to be
  interpreted uniformly across the base institution models in a Kripke
  model, or else by specific axioms in other cases (such as various
  properties of the Kripke frame, but not only). 
  The constraint axioms do not appear in the specifications as
  they are part of the definition of all Kripke models and therefore
  are common to all specifications; they are declared when defining
  the respective logic/institution. 

\end{itemize}

Let us see how these ideas are realised in the case of a concrete
example of a \Hspec\ specification.
The following \Hspec\ specification is that of a reconfigurable
calculator for natural numbers with a binary operation that in one state is sum
and in the other one is multiplication, an example which is discussed in
\cite{madeira-phd}.  

\begin{footnotesize}
\begin{verbatim}
spec Nat =
 logic : RigidCASL 
 rigid sort Nat
 rigid op 0 : Nat
 rigid op suc : Nat -> Nat
 op X : Nat * Nat -> Nat
end

spec Calc = 
  hlogic : HRigidCASLC 
  data Nat
  {
  nominals mult, sum
  modality shift : 2
  . mult \/ sum
  . @ sum
     : <shift> mult /\ [shift] mult
  . @ mult
     : <shift> sum /\ [shift] sum
  . @ mult : not sum
  . @ sum
     : forall m : Nat
         . X(m, 0) = m
  . @ sum
     : forall m, n : Nat
         . X(m, suc(n)) = suc (X(m, n))
  . @ mult
     : forallH m : Nat
      . forallH n : Nat 
      . existsH x : Nat
      . existsH y : Nat
         . X(m, n) = x /\ X(m, suc(n)) = y /\ <shift> X(x, m) = y 
 } 
end
\end{verbatim}
\end{footnotesize}
The first specification, at the level of the base institution,
declares the data of the natural numbers together with the binary
operation that will change modes that can be either interpreted as
addition or as multiplication.
This uses the CASL logic (essentially first order logic with partial
functions).  
Rigidity constraints are also specified at this stage in order to
prepare for the Kripke models in which the base models share their
underlying sets and some of the operations. 
The rigidity declarations do not have any semantic effect at the
level of the data (\texttt{Nat}), however they will have an effect at
the level of the hybridisation. 
This latter aspect, although does not introduce any error, is a little
``unclean''. 
It is not a kind of implementation shortcut, it rather comes from a
small gap in the theory. 
Currently, when building a hybridisation $\HI^C$, the signatures of the
base institution are preserved. 
Since rigidity of sorts and operations are in fact declarations at the
level of the base institution signatures, they have to be there
already in order to specify sharing constraints in the hybridisation.
A possible solution to this is to go more abstract about the
signatures of $\HI^C$ by specifying them abstractly together with a
projection functor to the signatures  of $\I$ (the base institution)
that may be subject to some axioms, in the style of how frame and
nominals extractions are defined in \cite{KripkeStrat}.   

The second specification is at the level of a hybridisation, which in
this case is \texttt{HRigidCASLC} (rigid sorts, rigid total functions
and the domain of each rigid partial function are interpreted
uniformly). 
Its definition does not appear in this specification as it resides in 
a library, being a predefined entity of \Hspec. 
However in this particular example we do not use any partial
functions. 
The data \texttt{Nat} is imported and nominals and the modality are 
declared. 
In the case of the modality note its arity 2. 
This is the part that declares the respective $\HI^C$ signature.
Then follows a series of axioms mainly regarding the dynamics of the
system.  
For example, the first axiom says that the Kripke frames have only two
elements.  
Note the two levels of quantifiers, \texttt{forall} is at the level
of the base institution while \texttt{forallH} and \texttt{existsH}
are quantifiers at the level of the hybridised institution. 
Because the base is a kind of fully fledged kind of first order logic,
its hybridisation $\HI^C$ differs substantially from the classical
first order hybrid logics. 

\subsection{Encoding into first-order logic}

The current version of \Hver\ contains only one method and tool that
is based upon a mathematical result that constitutes the main
achievement in \cite{EncHybrid}. 
That result represents an extension of the traditional translation 
of modal logic to first order logic \cite{bentham88} (for the hybrid
variant \cite{blackburn-seligman95}) to encodings of abstract 
hybridised institutions into first order logic.

That encoding uses the mathematical notion of \emph{comorphism}
\cite{jm-granada89,tarlecki95,tarlecki98,mossakowski96,goguen-rosu2000},
which is an important concept of institution theory. 
From the perspective of the mathematical structure, comorphisms are
just `homomorphisms of institutions'. 
So they are mappings between institutions that preserve the
mathematical structure of institutions. 

\

\noindent\fbox{%
    \parbox{\textwidth}{%
An \emph{institution comorphism}
$(\Phi, \alpha, \beta) \co \I \ra \I'$ consists of      
\begin{enumerate} 

\item a functor $\Phi \co \Sign \ra \Sign'$,

\item a natural transformation
  $\alpha \co \Sen \Rightarrow \Phi;\Sen'$, and

\item a natural transformation
  $\beta \co \Phi^\op ; \Mod' \Rightarrow \Mod$ 
\end{enumerate}

\noindent
such that the following \emph{satisfaction condition} holds 
\[
M' \models'_{\Phi(\Sigma)} \alpha_{\Sigma}(e) \mbox{ \ iff \ }
\beta_{\Sigma}(M') \models_{\Sigma}e
\]
for each signature $\Sigma\in |\Sign|$, for each $\Phi(\Sigma)$-model
$M'$, and each $\Sigma$-sentence $e$.
    }%
}

\

\noindent
While the $\alpha$ represents the translation of the syntax, the
$\beta$ represents the translation of the semantics. 
The \emph{satisfaction condition} ensures the mutual coherence between
these translations. 

Although comorphisms generally express an embedding relationship
between institutions,  they can also be used for `encoding' a `more
complex' institution $\I$ into a  `simpler' one $\I'$. 
In such encodings the structural complexity cost is shifted to the
mapping $\Phi$ on the signatures, thus $\Phi$ maps signatures of $\I$
to \emph{theories} of $\I'$ rather than signatures. 
This is why in the literature these are sometimes
\cite{mossakowski96,goguen-rosu2000}  called `theoroidal'
comorphisms. 
A \emph{theory} in $\I$ is just a specification in $\I$, i.e. a
signature $\Sigma$ plus a set $E$ of $\Sigma$-sentences. 
Technically speaking a `theoroidal' comorphism is in fact an ordinary
comorphism when we replace the institution $\I'$ with the
\emph{institution of its theories} $\I'^{\p}$.   
This is achieved through a general construction that can be applied to
absolutely any institution, in which the signatures of $\I'^{\p}$ are
the theories of $\I$. 
The details of this construction may be found in may places in the
literature, such as in \cite{iimt} (but under the name of the
institution of `presentations').  

Due to the generality of the construction of $\HI^C$, including its
parameters $\D$ and the constraint sub-functor on models, the
definition of the general encoding of $\HI^C$ into first order logic is
technically rather complex. 
Therefore we omit it here (for the details the interested reader has
to refer to \cite{EncHybrid}) and instead we present briefly its main
idea. 
The basis of the construction of the comorphism $\HI^C \to \FOL^\p$
(where $\FOL$ is the institution of first order logic in its many
sorted form) is a given encoding from the base institution to $\FOL$,
i.e. a comorphism $\I \to \FOL^\p$. 
This is considered abstractly, so it may vary, and in this way it
constitutes the main parameter of this construction. 
In practice these comorphisms may be well established translations. 
Then under some technical conditions -- mainly about quantifiers and
model constraints -- that are commonly satisfied in the applications, 
the comorphism $\I \to \FOL^\p$ is ``lifted'' to a comorphism  
$\HI^C \to \FOL^\p$. 

\subsubsection{How \Hver\ works in principle}

Suppose that we have a specification in \Hspec, which corresponds to a
theory $(\Delta,E)$ in some $\HI^C$. 
Suppose that we have a property $e$ that we have to check; this means
that we have to prove that $E \models_\Sigma e$ (which means that any
model of $E$ also satisfies $e$). 
In order to establish $E \models e$ we have to perform the following
steps: 
\begin{enumerate}

\item 
Then we translate $E$ and $e$ by using the comorphism $\HI^C \to
\FOL^\p$; this yields $\tilde{E}$ and $\tilde{e}$ in $\FOL$. 
Usually $\tilde{E}$ includes both the syntactic translation $\alpha(E)$
of $E$ and the sentences of the theory $\Phi(\Delta)$. 
Obviously in the case of $\tilde{e}$ it is not necessary to include
the latter sentences, so $\tilde{e} = \alpha(e)$.
 
\item 
If we had a theorem prover for $\FOL$ then this step would not be
necessary. 
But unfortunately, at least up to my knowledge, all major first order
logic theorem provers work with the \emph{unsorted} version of first
order logic.  So we translate again both $\tilde{E}$ and $\tilde{e}$
along a well known comorphism that encodes many sorted first order
logic into unsorted first order logic (the details of this comorphism
may be found for example in \cite{iimt}).  
We thus arrive at $\bar{E}$ and $\bar{e}$. 

\item
Now we employ a first order theorem prover and attempt to prove 
that $\bar{E} \models \bar{e}$. 

\item If we are successful with the previous task then we may conclude
  that $E \models e$.
  However this move backwards is not straightforward. 
  It holds as a consequence of an important property of comorphisms,
  namely that of \emph{conservativity}:

\

\noindent\fbox{%
    \parbox{\textwidth}{%
An institution comorphism
$(\Phi, \alpha, \beta) \co \I \ra \I'$ is \emph{conservative} when for
each $\Sigma$-model $M$ in $\I$, there exists a $\Phi(\Sigma)$-model
$M′$ in $\I′$ such that $M = \beta_\Sigma (M′)$. 
    }%
}

\

\noindent
In \cite{EncHybrid} we have shown that under some technical conditions
that are usually satisfied in the applications, the conservativity of
the comorphism $\HI^C \to \FOL^\p$ is inherited from the
conservativity of the base comorphism $\I \to \FOL^\p$. 
In order to complete the argument we still need that the encoding of
$\FOL$ to unsorted first order logic is conservative, which indeed
is. 

\end{enumerate}
In practice all the translations and the proofs are performed
automatically using tools. 
We will see more about this later on. 

\section{The current \Hlang\ implementation}

\subsection{ForVer and Hets}

In the year 2017 the author of this paper won funding for the project
proposal \emph{Formal Verification of Reconfigurable Systems}
(acronym: \emph{ForVer}) under a new funding scheme of a Romanian
government agency for funding of research that was dedicated to
experimental-demonstrative projects. 
That competition was highly competitive, with a succes rate of about
5\%, and the reviewers of the project proposals were selected from the
international scientific community. 
The goal of \emph{ForVer} was to realise the long term vision of
\Hlang\ and of the science behind it as a running prototype. 
Then the project hired Mihai Codescu for programming the prototype. 

The basic plan for this implementation of \Hlang\ was to rely on Hets
\cite{mossakowski-maeder-lutich2007}. 
Hets is a tool for heterogeneous multi-logic specification and
modeling of software systems and ontology development.  
In both these fields, there are a large number of logics and languages
in use, each better suited for a different task or providing a better
support for a different aspect of a complex system. 
Instead of trying to integrate the features of all these logics into a
single formalism, the paradigm of heterogenous multi-logic
specification is to integrate all logics by the means of a so-called
Grothendieck construction over a graph of logics and their
translations (captured as institutions and institution comorphisms,
respectively).    
Thus, for each logic we can make use of its dedicated syntax(es) and
proof tools. 
The specifier has the freedom to choose the logic that suits best the
problem to be solved, offers best tool support and according to the
degree of familiarity with a certain specification language. Hets
provides an implementation of this paradigm. 
Because of the multi-logic feature of \Hlang\ and also because of
\Hver\ is based logic encodings (translations), Hets appeared as 
suitable for a smooth implementation of a first prototype for \Hlang. 

\subsubsection{Grothendieck institutions} 

As mentioned above the foundation of Hets are the so-called
\emph{Grothendieck institutions}, which represents the ultimate
theoretical answer to the problem of heterogeneous multi-logic
specification. 
Instead of presenting the rather intricate technicalities of this
concept let us review how it was developed. 
This theory has been initially developed gradually within the context
of the CafeOBJ \cite{caferep}, which was the first
heterogeneous specification language.
A first attempt to address this heterogeneity institution
theoretically was in \cite{etm}. 
Then the late Professor Martin Hoffman, while writing a review for this
publication, suggested a construction on institutions similar to the
famous construction by Alexandre Grothendieck on categories
\cite{grothendieck63} originating from algebraic geometry. 
In the year 2000 this suggestion was realised in \cite{gi}, but that
was based on the original concept of homomorphism of institutions, the
\emph{institution morphisms} of \cite{ins}, which is somehow dual to
the concept of institution comorphisms that was discussed above. 
A few years later it was realised \cite{Mossakowski02b} that some
crucial proporties of Grothendieck institutions may be obtained more
smoothly by the same construction but based on institution comorphisms
rather than institution morphisms. 
The Hets system is based on the later version of Grothendieck
institutions. 

\subsection{The Hets implementation of \Hlang}

A more detailed description of this implementation, may be
found in \cite{codescu:hybrid-2019}. 
Note that this is an open implementation that supports further
enhancements. 
Here we review very briefly its main components:

\subsubsection{Syntactic support for the declarations of the parameters
  of the hybridisation process}

The considered parameters of the hybridisation are:
\begin{itemize}

\item The base institution. 
This is specified by using its internal Hets name, based on a Hets
qualification mechanism it is possible to select also a
sub-institution of an institution already implemented in Hets. 

\item The quantifier symbols.
These may be nominal symbols or classes of symbols that are specific
to the quantifications in the base institutions (such as constants,
rigid constants, total constants). 

\item The constraints on the Kripke models.
These are specified through a fixed grammar that cover two different
kinds of constraints: on the Kripke frames and on the interpretations
of symbols in the local models. 

\end{itemize}
The definitions of the hybridised logics are registered in such a way
that allows new additions of quantifications symbols and of
constraints. 

\subsubsection{Generic method for generating new instances of the Hets
  class \textbf{Logic}}

The generation of new instances of the Hets class logic \textbf{Logic}
is achieved on the basis of the definitions of hybridisations of
institutions by introducing Haskell polymorphic types for the
signatures of the base institution, for the nominals, and for the
modalities. 

\subsubsection{Support for structured specifications}

The specification structuring operators of \Hspec\ consist of unions,
specification translations alongs signature morphisms (which are
symbols renaming), and colimits of signatures.  
These are supported on the basis of a correspondence between the
structured \Hspec\ specifications and DOL \cite{DOL}, the structuring
language  supported in Hets.  
This correspondence is embedded in the concrete definition of \Hspec. 
 
\subsubsection{Support for \Hver}

There is a special declarative syntax for this that takes as
parameters the base theoroidal comorphism (from the base institution
$\I$ to $\FOL$) and the name of the respective hybridised logic
$\HI^C$.  
A generic method analises these definitions and generates Haskell code
containing a corresponding new instance of the type class for Hets
encodings.   
The compilation of this code makes the new encoding available for the
verification process, where the translation $\bar{E} \models \bar{e}$
of a goal $E \models e$ is passed to one of the first order logic
theorem provers of choice, such as SPASS \cite{spass3.5}, Vampire 
\cite{vampire}, E \cite{e1.8}. 

\subsubsection{Support for new logics}

In order to cover in \Hlang\ the institutions that are presented
as examples in \cite{EncHybrid} it was necessary to implement in Hets
the institution of the partial algebras with rigid symbols and first
order logic with rigid symbols (as a sub-institution of the former). 

\section{\Hlang\ at work}

A number of succesful case studies with \Hlang\ have been already
reported.  
In this section we present very briefly a couple of them.

\subsection{A steam boiler control}

The problem of this case study is a very notorious benchmark in
formal methods \cite{steam-boiler-book}. 
The case study with \Hlang\ has been reported in
\cite{codescu:hybrid-2019} and the \Hlang\ code is available at 
\begin{quotation}
\url{https://ontohub.org/forver/Sbcs.dol}
\end{quotation}

The \Hlang\ specification of the boiler control system illustrates
almost all of features of \Hlang. 
The base institution $\I$ is many-sorted first order logic, the
hybridised institution $\HI^C$ has quantifications over nominals and
rigid constants, and the constraints are given by the sharing of the 
domains and of the interpretations of rigid symbols. 

The properties that have been verified include changes of modes when
an event takes place and that in all states of the system the expected 
functionality takes place. 
In the \Hlang\ formalisation, the system has five modes (nominals) and
nine events (modalities). 

\subsection{A bike-sharing system design}

This case study has been reported in
\cite{Tutu-Chirita-Lopes-Fiadeiro:LSBSS-2019} and the \Hlang\ code is
available at 
\begin{quotation}
\url{https://ontohub.org/forver/BSS.dol}
\end{quotation}

It is based on a double hybridisation (hybridisation iterated twice)
the base level for the first level hybridisation being the (atomic
fragment) of propositional logic.  
\begin{itemize}

\item
The first level hybridisation has quantifications over nominals and
one constraint, namely that the interpretation of one of the
modalities (`parent') is  a forest (a set of disjoint trees). 

\item
The second level of hybridisation admits quantifications over first
level nominals (called ``actors'') as well as quantifications over
second level nominals (called ``configurations''). 
There is a sharing constraint: the first level Kripke models in a
second level Kripke model share the same underlying set of
``actors''. 

\end{itemize}
Since at the verification stage this modelling leads to some timeout
problems (due to a huge number of sentences obtained by the encoding
in $\FOL$), the first level of the hybridisation has been encoded in
an institution of relations. 

The first order theorem prover employed by this case study is SPASS.

\section{\Hlang\ in the future}

There are several directions that I see with respect to the future
evolution of \Hlang. 
\begin{itemize}

\item When conditions allow there should be a new implementation of
  \Hlang\ that is independent of Hets. 
  The reasons for this are manifold. 
  For example Hets is a rather big system and \Hlang\ relates only to
  a small part of Hets. 
  Such big systems are prone to errors that may easily affect the
  functionality of \Hlang. 
  Moreover \Hlang\ maintainers have little control on the evolution of
  Hets. 

\item \Hver\ should be enhanced with more tools and methods. 
  For example a direct tool based on proof systems at the level of
  hybridised institutions (so \emph{not} by translation) would be an
  welcome enhancement of \Hver.

\item Adding new base institutions and constraints to the current Hets
  implementation of \Hlang. 

\item Some slight upgrades of the foundations may be necessary in
  order to accomodate certain specification methodologies. 
  For instance, we have already discussed the issue of rigidity
  declarations at the level of the base institutions which may be
  solved by considering `projection' functors from the categories of
  the signatures of the hybridised institutions to the categories of
  the signatures of the base institution. 

\item More larger case studies should be developed with the aim to
  finally have \Hlang\ as an industrial tool. 

\end{itemize}

\bibliographystyle{apalike}
\bibliography{/Users/diacon/TEX/tex}%

\end{document}